\documentclass[12pt]{amsart}
\usepackage{amscd}
\usepackage[all]{xy}
\xyoption{dvips}
\newcommand{\PP}{{\mathcal{P}}}
\newcommand{\func}{{\mathcal{S}}}
\renewcommand{\SS}{{\operatorname{\tt SS}}}
\newcommand{\Set}{{\operatorname{\tt Set}}}
\newcommand{\Top}{{\operatorname{\tt Top}}}
\newcommand{\Hom}{{\operatorname{Hom}}}
\newcommand{\Map}{{\operatorname{Map}}}
\newcommand{\inject}{\hookrightarrow}
\newcommand{\obj}{{\operatorname{Obj}}}
\newcommand{\isom}{\cong}

\def\dlim_#1{\mathchoice{\underset{#1}{\varinjlim}}{\varinjlim_{#1}}{}{}}
\newcommand{\dtild}{{\tilde{\Delta}}}
\newcommand{\Dtild}{{\tilde{\nabla}}}
\newcommand{\Z}{{\mathbb{Z}}}
\newcommand{\R}{{\mathbb{R}}}
\newcommand{\PPP}{{\mathcal{PP}}}
\newcommand{\PPPo}{{\mathcal{PP}_1}}
\newcommand{\im}{\operatorname{Im}}

\newenvironment{Notation}%
  {\begin{list}{}%
      {\setlength{\leftmargin}{0pt}%
        \setlength{\labelwidth}{0pt}%
        \setlength{\itemsep}{3pt}%
      }%
  }%
{\end{list}}

\newtheorem{theorem}{Theorem}[section]
\newtheorem{lemma}[theorem]{Lemma}
\newtheorem{prop}[theorem]{Proposition}
\newtheorem{cor}[theorem]{Corollary}

\theoremstyle{definition}
\newtheorem{definition}[theorem]{Definition}
\theoremstyle{remark}
\newtheorem{remark}[theorem]{Remark}
\numberwithin{equation}{section}

\begin{document}
\title[Realization of simplicial and cyclic sets]{A simple approach to
  geometric realization\\of simplicial and cyclic sets}
\author{Amnon Besser}
\address{Department of Mathematics\\
Ben-Gurion University of the Negev\\
P.O.B. 653\\
Be'er-Sheva 84105\\
Israel
}
\maketitle

\section{Introduction}

The theory of simplicial sets and their realization is perhaps the
basis to the combinatorial approach to homotopy theory. The
construction of the geometric realization is very natural except for
one thing: It begins with the
postulation of the geometric realization of each of the ``standard
simplices'', to which no justification is given. It occurred to us that
by ``explaining'' the nature of this realization, the theory could
become even more natural. In particular, we wanted to find a
completely obvious proof to the fact that geometric realization
commutes with products in the right topology.

This explanation is achieved here
in the first two sections. We interpret the geometric realization of
the standard simplex $\Delta_n$ as
an appropriately topologization of the space of order preserving
maps from the 
unit interval to the ordered set $[n]:=\{1,\ldots, n\}$.

The definition
makes perfect sense for any finite partially ordered set (and in
fact more generally for categories). On the other hand, a partially
ordered set $P$ gives rise to a natural simplicial set and we show that
its realization is the same as that of $P$. On the point set level
this boils down to the fact that an order preserving map from the
unit interval to $P$ has to factor through a map $[n]\rightarrow P$
for some $[n]$. On the point set level it is again clear that for
partially ordered set the constructions of geometric realization and
of associated simplicial sets commute with products. Once we check
that this remains true in topology, the proof that geometric
realization of simplicial sets commutes with products is attained.

It turns out that a very similar idea works also for cyclic sets. One
needs to introduce a notion that plays the same role in the theory of
cyclic sets as that of partially ordered sets in the theory of
simplicial sets. This notion, introduced in section~\ref{sec:periodic}, is that
of a \emph{periodic partially ordered set}, which is just a partially
ordered set with an action of a free abelian group of finite rank, but
the definition of morphisms is a bit awkward. Among other things, we
show that the cyclic category of Connes is isomorphic to a certain
subcategory of periodic partially ordered set of degree (which is the
rank of the acting abelian group) $1$.

As the reader will notice, no face or degeneracy maps are mentioned
anywhere in the text (except just now). We find this to be one of the
more pleasing aspects of the theory.

We would like to thank Ed Efros and Igor Markov for their interest in
this work. In particular, the suggestion that the theory should extend
to cyclic sets was made by Igor Markov. We would like to thank UCLA,
the Max Planck Institute and the Isaac Newton Institute, where part of
this work was done. 

Note: This work was previously only available on my home page. It has
since been cited~\cite{Dri03} (who does a similar
construction, which as it turns out also appears in~\cite{Gra03})
and~\cite{Crh01}. In the interest of making it more easily accessible
I am submitting it to the archive without any changes. In particular,
I have not incorporated suggestions from a referee report that I
received for this paper. One important remark that the referee made
is that the alternative description of the cyclic category already
appeared in~\cite{Elm93}
 
\section{Preliminaries}
\label{sec:prelim}

In this section we briefly recall the theory of simplicial sets and their
realizations appearing in the standard literature.
\begin{Notation}
\item[$\PP$] the category of finite partially ordered sets with order
  preserving 
  maps as morphisms.
\item[{$[n]\in\PP$}] the ordered set $\{0<1<2<\cdots <n\}$,
for a non negative integer $n$.
\item[$I$] the closed unit interval $[0,1]$.
\item[$\Delta$] the subcategory of $\PP$ containing all the objects
  $[n]$.
\item[$\Set$] the category of sets.
\item[$\Top$] the category of topological spaces.
\item[$\SS$] the category of simplicial sets. By definition an object
  of $\SS$ is a contravariant functor $C:\Delta\rightarrow \Set$ and
  morphisms between objects are natural transformations.
\item[$\Delta_n$] (for $n\ge 0$) is the standard $n$-simplex.
  This is the simplicial set which is the contravariant functor on
  $\Delta$ represented by the object $[n]$. The map $[n]\rightarrow
  \Delta_n$ extends to a functor $\Delta \rightarrow \SS$ and one has
  $\Hom_{\Delta}([n],[m])\isom \Hom_{\SS}(\Delta_n,\Delta_m)$ (see
  below lemma~\ref{partially}).
\item[$|\Delta_n|$] The topological standard $n$-simplex, defined as
  \begin{equation*}
    |\Delta_n|:=\{(x_0,x_1,\ldots,x_n)\in \mathbb{R}^{n+1},\
  x_i\ge 0 \text{ and } \sum x_i=1\}.
  \end{equation*}
  with the subspace topology. The map $[n]\rightarrow |\Delta_n|$
  extends to a functor $\Delta \rightarrow \Top$
  (see~\cite[p. 3]{Gor-Jar97}): For a map $\theta:[n]\rightarrow [m]$
  corresponds a map $\theta_\ast:|\Delta_n|\rightarrow |\Delta_m|$
  defined by $\theta_\ast(t_0,\ldots,t_n)=(s_0,\ldots, s_m)$ with
  \begin{equation*}
    s_i=\sum_{j=\theta^{-1}(i)} t_j.
  \end{equation*}
\item[$|\cdot|:\SS\rightarrow\Top$] The functor of geometric
  realization. We will use the following convenient
  definition~\cite{Jar89,Gor-Jar97}:
  \[ 
    |C|=\dlim_{\Delta_n \rightarrow C} |\Delta_n|.
  \]
  The limit is taken over the category $\Delta\downarrow C$
  whose objects are maps
  of simplicial sets $f:\Delta_n \rightarrow C$ and whose
  morphisms are commuting triangles
  \begin{equation*}
    \xymatrix{
      \Delta_n \ar[rr]\ar[dr]^{f} & &\Delta_m\ar[dl]_{g}\\
      & C & 
    }
  \end{equation*}
\end{Notation}

The categories $\PP$ and $\SS$ have products - the product of $P,\ Q\in \PP$ is the set
$P\times Q$ with the order $(x,y)\ge (z,w)$ if $x\ge z$ and $y\ge w$.
The product of the simplicial sets $C$ and $D$ is the
functor given on objects by $C\times D([n])= C([n])\times D([n])$.

\section{Realization of simplicial sets}
\label{sec:realization}

\begin{definition}
  The \emph{simplicial realization} of a partially ordered set $P$ is the
  simplicial set $\func(P)$ defined by 
  \[
    \func(P)([n])=\Hom_\PP([n],P).
  \]
\end{definition}
The assignment $P\rightarrow \func(P)$ defines a functor
$\func:\PP \rightarrow \SS$.
Clearly $\func([n])=\Delta_n$.
It is also clear that 
\begin{equation}\label{funcproduct}
\func(P\times Q)=\func(P)\times \func(Q).
\end{equation}
\begin{lemma}\label{partially}
  The natural map
  \begin{equation}\label{natural-map}
    \Hom_\PP(P,Q) \rightarrow  \Hom_\SS (\func(P),\func(Q))
  \end{equation}
        is a bijection.
\end{lemma}
\begin{proof}
  A map $F:\func(P) \rightarrow \func(Q)$ associates to every order
  preserving map $f:[n]\rightarrow P$ an order preserving map
  $g=F(f)$, $g:[n]\rightarrow Q$ in such a way that for
  $\theta:[m]\rightarrow [n]$ one has $F(f\circ \theta)=F(f)\circ
  \theta$. By considering the case $n=0$ we see that $F$ induces a
  unique map $\tilde{F}:P\rightarrow Q$ such that $F(f)= f\circ
  \tilde{F}$ for any map $f:[0]\rightarrow P$. By considering
  $\theta_i: [0]\rightarrow [n]$ sending $0$ to $i\in [n]$ we see that
  $g(i)=\tilde{F}(f(i))$ for all $i$, i.e., that $F(f)=\tilde{F}\circ
  f$. For $F$ to send an order preserving $f$ to an order preserving
  $F(f)$ it is necessary and sufficient (consider the case $n=1$) that
  $\tilde{F}$ is order preserving. It is clear that $F\rightarrow
  \tilde{F}$ is an inverse to \eqref{natural-map}.
\end{proof}
One can define the geometric
realization of a partially ordered set $P$ to be the geometric
realization of $\func(P)$. We will later show that this is equivalent
to the following direct definition:
\begin{definition}\label{geomreal}
  Let $P$ be a finite partially ordered set. Consider $P$ as a metric space by
  the standard discrete metric $d(x,y)=1-\delta_{xy}$ for any $x$ and
  $y$ in $P$. The {\em geometric realization\/} of $P$, denoted $|P|$,
  is the set 
  of all order preserving upper semi-continuous maps
  $f:I\rightarrow P$. We give $|P|$ the
  structure of a metric space by the metric
  \[
  d(f,g)=\int_0^1 d(f(t),g(t))dt.
  \]
\end{definition}
One checks immediately that $|\ |$ is a functor from $\PP$ to $\Top$.
\begin{lemma}\label{mainlemma}
  For any $n\ge 0$ the geometric realization of $[n]$ is homeomorphic
  to the geometric 
  realization of $\Delta_n$. The homeomorphism is given by the map
  \[ f\rightarrow
  (\mu(f^{-1}(0)),\mu(f^{-1}(1)),\ldots,\mu(f^{-1}(n)))\] with $\mu$
  the standard measure on $I$. Moreover, this homeomorphism is
  compatible with the maps induced by morphisms in $\Delta$: if
  $\theta:[n] \rightarrow[m]$ is a morphism then we have the
  commutative diagram
  \begin{equation*}
    \begin{CD}
      |[n]| @>>> |\Delta_n| \\
      @V{|\theta|}VV @VV{\theta_\ast}V \\
      |[m]| @>>> |\Delta_m|
    \end{CD}
  \end{equation*}
\end{lemma}
\begin{proof}
  To check continuity
  it is enough to note that $\mu(f^{-1}(j)) =1-d(f,j)$ where here the
  last $j$ means the constant function $j$ on $I$.
  The inverse map is given by
  $\underline{t}=(t_0,\ldots,t_n)\rightarrow f_{\underline{t}}$, with
  $f_{\underline{t}}(t)=j$ if $t_0+\cdots +t_j  \le t <  t_0+\cdots
  +t_{j+1}$. To see that this inverse map is continuous as well, suppose we have
  $\underline{t}$ and $\underline{t^\prime}$ with
  $|t_i-t_i^\prime |< \epsilon$. Letting $r_j=t_0+\cdots +t_j$ and
  similarly for $\underline{t^\prime}$, it follows that 
  \begin{equation*}
    \mu([r_j,r_{j+1}]-[r_j^\prime, r_{j+1}^\prime])\le 2(n+1)\epsilon,
  \end{equation*}
  and since when $t\in [r_j,r_{j+1}]$, $f_{\underline{t}}$ differs from
  $f_{\underline{t}^\prime}$ only in the above set, we find that
  $d(f_{\underline{t}},f_{\underline{t}^\prime})\le 2(n+1)^2
  \epsilon$. The compatibility with the morphisms in $\Delta$ is
  easily established.
\end{proof}
\begin{prop}\label{peqsp}
  There is a functorial isomorphism, for any $P\in \PP$:
  \[
    |P|\isom |\func(P)|.
  \]
\end{prop}
\begin{proof}
  For any $[n]\rightarrow P$, we obtain by functoriality a map
  $|[n]|\rightarrow |P|$ (we remark that this map takes $(I\rightarrow
  [n])\in |[n]|$ to the composition $I\rightarrow[n]\rightarrow
  P$). Therefore, there is a map in $\Top$ 
  \begin{equation}\label{natmap}
    \dlim_{ [n]\rightarrow P }|[n]|\rightarrow |P|,
  \end{equation}
  where the limit is over the category $\Delta\downarrow P$ whose
  object are maps 
  $[n]\rightarrow P$ and whose morphisms are the obvious triangles. It
  follows from Yoneda's lemma (or from lemma~\ref{partially}) that this
  last category is isomorphic to
  $\Delta\downarrow \func(P)$. We will show that the map \eqref{natmap} is a
  homeomorphism and will therefore be done by
  lemma~\ref{mainlemma}. To construct an inverse to \eqref{natmap} we make
  the key observation that any order preserving map
  $f:I\rightarrow P$ factors as $I\xrightarrow{f_0} [n]
  \xrightarrow{f_1} P$ for 
  some $n$. In fact, we may choose a canonical such decomposition: The
  image of $f$ in $P$ is a finite totally ordered set. There is therefore a
  unique $n$ and a unique order preserving bijection $[n]\isom \im(f)$ which
  allows us to factor $f$ through $[n]$. Let us write this decomposition, by
  abuse of notation, as $I\xrightarrow{f} \im(f) \inject P$. In this way we
  obtain from $f$ an element of $|\im(f)|$ and therefore of $\dlim_{
    [n]\rightarrow P }|[n]|$. This gives a well defined map of sets
  $|P| \rightarrow \dlim_{[n]\rightarrow P }|[n]|$. It is clear that the
  composite map  $P \rightarrow \dlim_{[n]\rightarrow P }|[n]|
  \rightarrow |P|$ is the identity map. We show that one gets the
  identity in the reverse 
  direction as well: Start with a sequence of maps $I\xrightarrow{u}
  [n] \xrightarrow{v} P$ representing an object $[n]\rightarrow P$ of
  $\Delta \downarrow P$ and an element of $|[n|$. If we map this
  element to $\dlim_{[n]\rightarrow P }|[n]|$, then to $|P|$ and then
  back to  $\dlim_{[n]\rightarrow P }|[n]|$, we obtain the image in
  $\dlim_{[n]\rightarrow P }|[n]|$ of the sequence $I\rightarrow \im
  v\circ u \inject P$ and we need to show that this has the same image
  as the original sequence. This is done by observing that both
  sequences have the same image in the realization of $\im v \inject
  P$.
  To show that the bijection \eqref{natmap} is a homeomorphism, it is
  enough to show that
   $\dlim_{ [n]\rightarrow P  }|[n]|$ is compact. This may be achieved
  by noticing that the subcategory of $\Delta\downarrow P$ consisting
  of injective maps $[n]\rightarrow P$ is cofinal. This implies by
  \cite[Theorem IX.3.1]{Mac71} that one may take the limit only over this
  subcategory, which 
  has only a finite number of objects.

\end{proof}
\begin{lemma}\label{lastlab}
  For any $P,\ Q\in \PP$ one has 
  \begin{equation}\label{geoproduct}
  |P\times Q|\isom |P|\times |Q|.
  \end{equation}
\end{lemma}
\begin{proof}
  The projections $P\times Q \rightarrow P$ and $P\times Q \rightarrow Q$
  induce continuous maps
  $|P\times Q|\rightarrow |P|$ and  $|P\times Q|\rightarrow |Q|$
  and therefore a
  continuous map  $|P\times Q|\rightarrow |P|\times |Q|$. On the level
  of sets it is clear that this map is a
  bijection. Since the two sides are compact by the proof of the
  previous proposition, we are done.
\end{proof}
\begin{prop}
  If $P,\ Q\in \PP$ then $|\func(P)\times \func(Q)|\isom |\func(P)|\times
  |\func(Q)|$. 
\end{prop}
\begin{proof}
  This follows from proposition~\ref{peqsp}, lemma~\ref{lastlab} and equation
  \eqref{funcproduct}.
\end{proof}
\begin{cor}\label{maincor}
  For any non negative integers $n$ and $m$ we have $|\Delta_n\times
  \Delta_m| \isom |\Delta_n|\times  |\Delta_m|$
\end{cor}
\begin{proof}
  This is just the last proposition when choosing $[n]$ and $[m]$ for
  $P$ and $Q$.
\end{proof}
\begin{remark}
  It is perhaps instructive to ``see'' the corollary in the case
  $n=m=1$. In this case $P$ and $Q$ are both equal to the ordered set
  $\{0<1\}$ and $P\times Q$ is therefore the poset
  \begin{equation*}
  \xymatrix@ur{
    (0,0) \ar@{}[r]|*[@]{<} \ar@{}[d]|*[@]{<}& (0,1) \ar@{}[d]|*[@]{<}\\
    (1,0) \ar@{}[r]|*[@]{<}& (1,1)
    }
  \end{equation*}
  The image of an order preserving map $I\rightarrow P\times Q$ can
  either be the ``lower path'' $\{(0,0)<(1,0)<(1,1)\}$, with boundary
  cases where the image is a subset of the path, or the ``upper path''
  $\{(0,0)<(0,1)<(1,1)\}$ and its subsets. Each of these paths corresponds
  to a triangle $|[2]|$ and the two are glued by the common edge which
  corresponds to the common subset $\{(0,0)<(1,1)\}$. The realization
  is therefore a square.
\end{remark}
\begin{remark}\label{product}
  corollary \ref{maincor} implies by a standard argument (see the proof of 
\cite[Theorem III.3.1]{Gab-Zis67}) that for any two simplicial sets $C$ and $D$
 one has
 $|C\times D|\isom |C|\times |D|$ provided one works in the category of
  Kelly spaces.
\end{remark}

\section{Periodic partially ordered sets}
\label{sec:periodic}

In this section we introduce a notion which is the cyclic equivalent of a
partially ordered set. This may look a bit unnatural. The usefulness
of the definition will be apparent in following sections.
\begin{definition}
  A \emph{periodic partially ordered set} (ppset for short), of degree $k>0$,
  is a partially 
  ordered set $P$ together with an action of $\Z^k$ by order
  preserving transformations. In
  other words, to give a ppset of degree $k$ is to give the poset $P$
  together with a 
  collection of commuting order preserving automorphisms
  $T_1,\ldots,T_k$ (We will call the automorphisms
  obtained by the action of $\Z^k$ the {\em shifts\/} of $P$). A
  \emph{map} $f:P\rightarrow Q$ between
  two ppsets is a function $f$ which is order
  preserving and satisfies the relation
  \[
    f\circ (1,1,\ldots,1)=(1,1,\ldots,1)\circ f
  \]
  where the symbol $(1,1,\ldots,1)$ represents the action of the
  appropriate group element (which could of course be different in $P$
  and $Q$). We write $\Map(P,Q)$ for the set of maps between ppsets
  $P$ and $Q$. A \emph{morphism} between two ppsets $P$ and $Q$ is an
  equivalence class of maps. The equivalence  is given by precomposition
  with the shifts of $P$ and postcomposition with the shifts of
  $Q$. The set of morphisms between $P$ and $Q$ is denoted $\Hom(P,Q)$.
\end{definition}
Unfortunately the collection of ppsets is \emph{not} a category with
respect to morphisms, because it is easy to see that composition does
not preserve the equivalence relation in general. On the other hand,
the subcollection of ppsets of degree $1$ does become a category this
way, because the commutation with $T_1$ implies that it is enough to
divide by the shifts on one side.
\begin{definition}
  The category $\PPPo$ is the category whose objects are degree $1$
  ppsets  and whose morphisms are morphisms of ppsets.
\end{definition}
The collection of all ppsets may be given the following formal
structure:
\begin{definition}
  Let $C$ be a category. A \emph{module} $M$ over $C$ is made up of
  the following data:
  \begin{enumerate}
  \item A class of objects $\obj(M)$.
  \item For any $P\in C$ and $Q\in M$, a set $\Hom(P,Q)$.
  \item For any $P,R\in C$ and $Q\in M$ a composition map
    $\Hom(R,P)\times \Hom(P,Q)\rightarrow \Hom(R,Q)$.
  \end{enumerate}
  The composition of morphisms should satisfy all the axioms
  of category theory when applicable. If $M^\prime$ is a $C^\prime$
  module then a 
  functor $F$ from $M$ to $M^\prime$ consists of
  \begin{enumerate}
  \item A functor $F:C\rightarrow C^\prime$.
  \item a map $F:\obj(M)\rightarrow \obj(M^\prime)$.
  \item For any $P\in C$, $Q\in M$, a map $F:\Hom(P,Q)\rightarrow
    \Hom(F(P), F(Q))$ satisfying the obvious compatibilities.
  \end{enumerate}
\end{definition}
The simplest example of a module occurs when $M$ is a category and $C$ is a
subcategory. A structure of a $C$ module on $M$ is obtained by forgetting
the morphisms between objects of $M$ and remembering only the ones
between objects of $C$ and objects of $M$. The following proposition,
whose proof is left to the reader,
is the reason for the introduction of the notion of a module:
\begin{prop}
  The collection $\PPP$ of ppsets of arbitrary degree forms a module
  over the category $\PPPo$. For $P\in \PPPo$ and $Q\in \PPP$,
  the set $\Hom(P,Q)$ is that of ppset morphisms between $P$ and $Q$
  and the composition is induced by the composition of maps of ppsets.
\end{prop}
\begin{definition}
  For $n\ge 0$ the \emph{standard ppset} $[[n]]$ is the degree $1$
  ppset defined as follows: As a partially
  ordered set it is $\Z$ and the generator
  $T_1$ is translation by $n+1$. 
\end{definition}
\begin{definition}\label{newdef}
  The category $\Dtild$ is the subcategory of $\PPPo$ whose objects are
  $[[n]]$ for $n=0,1,\ldots$.
\end{definition}

Like the case of posets, we need a notion of a product of
ppsets. Although the collection of ppsets is not a category, a notion
of products still exists as follows:
\begin{definition}
  Suppose $C$ is a category and $M$ is a $C$-module.
  \begin{enumerate}
  \item When $Q\in M$, the functor $C\rightarrow \Set$ given on
    objects by $P\rightarrow \Hom(P,Q)$ (visibly seen to be a
    functor), is called the functor represented by $Q$.
  \item An object $Q_3\in M$ is \emph{a} product of objects $Q_1$ and $Q_2$
    in $M$ if $Q_3$ represents the product of the functors represented
    by $Q_1$ and $Q_2$ (note that, unlike the case of categories, the
    product need not be unique).
  \end{enumerate}
\end{definition}
\begin{prop}
  The $\PPPo$ module $\PPP$ has products. A product of $P$ and a $\Z^k$
  action with $Q$ and a $\Z^m$ action is given by the product poset $P\times
  Q$ together with the product action of $\Z^k\times \Z^m\isom
  Z^{k+m}$.
\end{prop}
\begin{proof}
  It is immediate to see that, if $R\in \PPP$, then $\Map(R,P\times
  Q)$ is naturally isomorphic to $\Map(R,P)\times \Map(R,Q)$. This
  isomorphism descends to morphisms if $R\in \PPPo$ as
  \begin{equation*}
    \Hom(R,P\times Q)=\Map(R,P\times Q)/\Z^{k+m}\isom \Map(R,P)/\Z^k \times
    \Map(R,Q)/\Z^m. 
  \end{equation*}
  Note that the crucial point here is that we only need to divide on
  the right side.
\end{proof}

\section{cyclic sets}

In this section we want to show that the category $\Dtild$ introduced
in the last section is in fact isomorphic to the cyclic category of
Connes. We will later use this to provide a new definition of the
geometric realization of cyclic sets.
We will take~\cite{Jon87} as our basic reference on the cyclic
category and cyclic sets. A cyclic
set is a contravariant functor $F:\dtild\rightarrow \Set$, where
the category $\dtild$ will be defined below. For each $[n]\in \Delta$
let $K([n])$ be the cyclic group $\Z/(n+1)$ thought of as the group of
cyclic permutations of $[n]$.
\begin{definition}\label{dtilddef}
  The category $\dtild$ has the same objects as $\Delta$. The sets of
  morphisms are given by
  \begin{equation*}
    \Hom_{\dtild}([n],[m])=\Hom_\Delta([n],[m])\times K([n]),
  \end{equation*}
  and the composition of $(\phi,u)\in \Hom_\Delta([m],[k])\times
  K([m])$ and $(\chi,v)\in \Hom_\Delta([n],[m])\times K([n])$  is defined by
  \begin{equation*}
    (\phi,u)\circ (\chi,v)=(\phi\circ u_\ast \chi,\chi^\ast u\circ v).
  \end{equation*}
  Here, the order preserving map $u_\ast \chi$ and the cyclic
  permutation $\chi^\ast u$ are defined as follows: For $i\in [m]$ let
  $A_i=\chi^{-1}(i)$ and let $B_{u(i)}=A_i$. Then $[n]$ can be given a
  new ordering as the ordered union of the posets $B_i$. The cyclic
  permutation 
  $\chi^\ast u$ is the unique automorphism of $[n]$ which is order
  preserving from $[n]$ with the new ordering to $[n]$ with the standard
  one. Finally  $\chi^\ast u:=u\chi(\chi^\ast u)^{-1}$.
\end{definition}

The following notation and trivial remarks will be useful when we
compare this definition with an alternative one: For $\chi\in
\Hom_\Delta([n],[m])$ and $v\in K([n])$ we write $\chi, v\in
\Hom_{\dtild}([n],[m])$ for $(\chi,\operatorname{id})$ and
$(\operatorname{id},v)$ respectively. Then
we immediately see that $(\chi, v)=\chi \circ v$, that composition in
$\dtild$ of morphisms of $\Delta$ (resp. elements of $K([n])$) is the
same as the standard composition and that the composition rule is
uniquely determined by the relation
\begin{equation*}
  u\circ \chi=u_\ast \chi \circ \chi^\ast u.
\end{equation*}
We also wish to recall for future use the following
observation
\begin{lemma}\label{compinterp}
  With the notation as in definition~\textup{\ref{dtilddef}},
  $\chi^\ast u$ is given by
  translation by $-k$, where $k=\min B_i$ and $i$ is the smallest
  index for which  $B_i$ is not empty.
\end{lemma}
\begin{proof}
  Indeed, this $k$ becomes the smallest element in $[n]$ with the new
  ordering and is therefore mapped to $0$ by  $\chi^\ast u$.
\end{proof}

\begin{theorem}\label{catequiv}
  The categories $\Dtild$ and $\dtild$ are isomorphic. The isomorphism
  is given  by the functor $F:\dtild\rightarrow \Dtild$ defined as
  follows: On objects $F([n])=[[n]]$. For $[n]$ and $[m]$ in $\dtild$
  and $\chi \in \Hom_{\Delta}([n],[m])\subset \Hom_{\dtild}([n],[m])$,
  $F(\chi)$ is defined by
  \begin{equation*}
    F(\chi)(a+r\cdot (n+1))=\chi(a)+r\cdot (m+1),\qquad 0\le
    a<n+1,r\in \Z.
  \end{equation*}
  For $u\in K([n])$, $F(u):[[n]]\rightarrow [[n]]$ is translation by
  $u$. In general, $F(\chi\circ u)= F(\chi)\circ F(u)$.
\end{theorem}
\begin{proof}
  Using the remarks after definition~\ref{dtilddef}, one easily checks that
  $F$ is well defined (as a map on morphisms without claiming anything
  about multiplicativity). We begin the proof by constructing an inverse
  $G$ to $F$. Clearly on objects we must put $G([[n]])=[n]$. On
  morphisms $G$ is defined as follows: Let $f\in
  \Hom_{\Dtild}([[n]],[[m]])$, represented by $f:\Z\rightarrow\Z$ such that $f$
  is order preserving and $f(x+n+1)=f(x)+m+1$ for all $x\in \Z$. We
  set $G(f)=(G_1(f),G_2(f))=(\chi,v)$, where $v=-\inf
  f^{-1}\{0,1,2,\ldots\})$ and 
  $\chi(x)=f(x-v)$ for $x\in [n]$. If we change $f$ to an equivalent
  map $f^\prime(x)=f(x)+r\cdot (m+1)=f(x+ r\cdot(n+1))$ with $r\in
  \Z$, then $G_2(f^\prime)=G_2(f)- r\cdot (n+1)$ and therefore
  $G_1(f)=G_2(f)$. It is also easy to see that $\chi$ maps $[n]$ into
  $[m]$. Indeed, by definition $f(-v)\ge 0$. If $f(n-v)\ge m+1$, then
  $f(-v-1) = f(n-v)-m-1\ge 0$, contradicting the minimality of
  $-v$. Since $f$ is order preserving we have for $i\in [n]$, $0\le
  f(i-v) < m+1$.
  Therefore, $G$ is well defined.

  It is straightforward to check that $G$ is inverse to $F$. It is
  also clear that $G$ and $F$ respect 
  composition on morphisms in $\Delta$ and in $K([n])$.
  To complete the proof, it is enough to show 
  that 
  \begin{equation}\label{commutiden}
    G(F(u)\circ F(\chi))= u_\ast \chi\circ \chi^\ast u.
  \end{equation}
  Indeed, If this is the case, then since $F$ and $G$ are inverse we
  will obtain 
  \begin{equation*}
    F(u)\circ F(\chi)=F(u_\ast \chi\circ \chi^\ast u)=F(u_\ast
    \chi)\circ F(\chi^\ast u),
  \end{equation*}
 and in general we will have
  \begin{equation*}
    \begin{split}
    F((\phi\circ u)\circ (\chi \circ v))=F((\phi \circ u_\ast
    \chi)\circ (\chi^\ast u\circ v))=F(\phi)\circ [F(u_\ast \chi)\circ
    F(\chi^\ast u)]\circ F(v)\\=F(\phi)\circ [F(u)\circ F(\chi)]\circ
    F(v)=F(\phi\circ u)\circ F(\chi\circ v).
    \end{split}
  \end{equation*}
  Because of the way $u_\ast \chi$ is defined in~\ref{dtilddef} it is enough
  to check the identity \eqref{commutiden} for the $K([n])$ part, i.e., to show
  that $G_2( F(u)\circ F(\chi))=\chi^\ast u$.
  Set $f=F(u)\circ F(\chi)$. It is clear that $f(\Z)\cap [m]\ne \emptyset$
  hence $G_2(f)= -\inf (f^{-1}[m])$. It is now easy to check that
  for $i\in [m]$

  \begin{equation*}
    f^{-1}(i)=
    \begin{cases}
      \chi^{-1}(u^{-1}(i))&\text{ if } i\ge u;\\
      \chi^{-1}(u^{-1}(i))-n-1 &\text{ if } i<u.
    \end{cases}
  \end{equation*}
  From this one sees that, with the notation of definition~\ref{dtilddef},
  $\chi^{-1}(u^{-1}(i))=B_i$. Now one readily sees from the fact that
  $f$ is order preserving that the minimal
  value of $f^{-1}([m])$ is, modulo $n+1$, the smallest element in the
  $B_i$ with the smallest index $i$ for which $B_i$ is not
  empty. The theorem now follows from lemma~\ref{compinterp}
\end{proof}
\begin{remark}
  With the equivalent definition~\ref{newdef} of the cyclic category,
  it is very easy to see 
  that it is self dual: Given $[[n]]$ in $\Dtild$, consider the set
  $\operatorname{Map}([[n]],[[0]])$. An element $f\in
  \Map([[n]],[[0]])$ satisfies $f(x+n+1)=f(x)+1$ and is therefore
  surjective and uniquely determined by the number 
  $i=\inf f^{-1}(0)$. Clearly, the set
  $\operatorname{Map}([[n]],[[0]])$ is totally ordered by the relation
  of inequality of functions and has a shift given by pre or
  post composition with shifts as usual, making it an object of
  $\PPPo$. The unique $f$ corresponding 
  to $i$ is the map 
  $f_i(x)= [(x-i)/(n+1)]$ (here $[~]$ denotes the integer part
  function). This gives a bijection
  $\operatorname{Map}([[n]],[[0]])\rightarrow \Z$ which is even order
  preserving. The shift corresponds to addition of
  $n+1$. Thus, as an object of $\PPPo$,
  $\operatorname{Map}([[n]],[[0]])\isom[[n]]$. One can easily check
  that this makes $[[n]]\rightarrow 
  \operatorname{Map}([[n]],[[0]])$ into a contravariant isomorphism
  from $\Dtild$ to itself.
\end{remark}

\section{Realization of cyclic sets}
\label{sec:realcyc}

In this section we wish to mimic the constructions and results of
section~\ref{sec:realization} for cyclic sets.
We interpret the unit circle $S^1$ as an object in $\PPP$. The
underlying ordered set is $\R$ - the set of real numbers. The shift
operator is translation by $1$.

The analogue of a finite poset is the notion of a compact ppset while that
of a totally ordered set is the notion of an archimedean ppset.
\begin{definition}
  a ppset $P$ of degree $k$ is compact if $\Z^k\backslash P$ is finite.
\end{definition}
\begin{definition}
  A ppset $P$ is called archimedean if it is of degree $1$, totally
  ordered and for any 
  $x$ and $y$ in $P$ 
  there is some $n\in \Z$ such that $T_1^n x > y$.
\end{definition}
\begin{lemma}
  If $P$ is an archimedean ppset, then either for all $x\in P$,
  $T_1(x)>x$, or for all $x\in P$, $T_1(x)<x$.
\end{lemma}
\begin{proof}
  Notice that an archimedean ppset must have more than one element.
  It is enough to prove that it is impossible to have $x\ne y\in P$ such
  that $T_1(x)\ge x$ while $T_1(y)\le y$. Suppose that we have these
  $x$ and $y$. Clearly $T_1^n(x)\ge x$ and $y\ge T_1^n(y)$ for $n>0$
  with the reverse inequalities holding for $n<0$. Suppose
  $y>x$. Since $P$ is archimedean, we have some $n$ such that $T_1^n
  y<x$. We must have $n>0$. But then $x\ge T_1^{-n}x > y$ and we arrive
  at a contradiction. Similarly, if $x>y$ we have, for some $n<0$,
  $T_1^n x< y$ which implies that $x<T_1^{-n}y \le y$.
\end{proof}
\begin{definition}
  We will say that an archimedean ppset $P$ is positive
  (resp. negative) if for one (hence any) $x\in P$, $T_1(x)> x$
  (resp. $T_1(x)< x$).
\end{definition}
\begin{definition}
  The geometric realization of a compact $P\in \PPP$ is the topological space
  $||P||$ defined as followed: As a set $||P||$ is the subset of
  $\Hom_\PPP(S^1,P)$ of morphisms where the underlying map of sets is
  upper semicontinuous when $P$ is taken to have the discrete topology.
  restriction of maps to $I\subset \R$ and projection on $\Z^k\backslash P$
  gives an embedding of $||P||$ into a 
  space of maps from $I$ to $\Z^k\backslash P$ and the metric on $||P||$ is
  then defined similarly to the one in definition~\ref{geomreal}.
\end{definition}
The precomposition with translation on $\R$ makes $||P||$ naturally
into a space with a circle action and it is immediate that $||\cdot||$
is a functor from $\PPP$ to the category of topological spaces with
circle action, considered as a module over itself.
\begin{prop}
  The geometric realization of $[[n]]$ is homeomorphic to the product
  $|\Delta_n|\times S^1 $ where the circle action fixes $|\Delta_n|$
  and acts in the obvious way on $S^1$.
\end{prop}
\begin{proof}
  This is proved in much the same way as theorem~\ref{catequiv}: To an order
  preserving map
  $f:\R \rightarrow \Z$ satisfying $f(x+1)=f(x)+n+1$, which represents
  an element of $||[[n]]||$, we associate the pair $(\phi,s)\in
  |[n]|\times S^1 $ as follows: $s=-\inf f^{-1}(\{0,1,2,\ldots\})$ and
  $\phi$ is given by $\phi(x)=f(x-s)$. The inverse map is given by
  $(\phi,s)\rightarrow f(x) = \phi(\{x+s\})+[x+s]$ where $\{~\}$ and
  $[~]$ denote fractional and integer values respectively. One only
  needs to check that these maps are continuous, which is easy.
\end{proof}

The theory can now be developed exactly as in the simplicial case. We
only sketch the proofs as needed.

\begin{definition}
  The cyclic realization $\mathcal{C}(P)$ of a ppset $P$ is the
  restriction to $\Dtild$ of 
  the functor on $\PPPo$ represented by $P$.
\end{definition}
\begin{definition}
  The standard cyclic $n$ simplex is $\dtild_n:=\mathcal{C}([[n]])$. 
\end{definition}
\begin{definition}
  The geometric realization of a cyclic set $C$ is given by:
  \begin{equation*}
    ||C||=\dlim_{\dtild_n \rightarrow C} ||[[n]]||.
  \end{equation*}
  The limit is taken in the category of topological spaces with circle
  action.
\end{definition}
The analogue of proposition~\ref{peqsp} is now
\begin{prop}
  For any compact ppset $P$ we have $||\mathcal{C}(P)||\isom ||P||$.
\end{prop}
\begin{proof}
  The proof is essentially the same as that of proposition
  \ref{peqsp}. The only point 
  which is maybe not obvious is the analogue of the fact that the
  image of $I$ by an order preserving map is isomorphic to $[n]$ for
  some $n$. This is provided by lemma \ref{analog} below, together with the
  obvious fact that the image of $S^1$ under a ppset map is
  positive archimedean and is compact if the target is.
\end{proof}
\begin{lemma}\label{analog}
  A positive archimedean compact ppset $P$ is isomorphic to $[[n]]$
  where $n+1$ is 
  the cardinality of $\Z\backslash P$.
\end{lemma}
\begin{proof}
  Let $\pi:P \rightarrow\Z\backslash P$ be the projection. Choose an
  element in $P$ and call it $0$. We can construct a section
  $s:\Z\backslash P \rightarrow P$ to $\pi$ in the following way: to
  each $\tilde{x}\in \Z\backslash P$, the set $A_{\tilde{x}}:=\{x\in
  \pi^{-1}(\tilde{x}),x\ge 0\}$ has a smallest element. To see this choose
  $x_0\in A_{\tilde{x}}$.
  Then $A_{\tilde{x}}=\{T_1^n(x_0)\ge 0\}$ and the set of $n$ for which
  $T_1^n(x_0)\ge 0$
  is bounded from below. Let $s(\tilde{x})$ be that smallest
  element. The image set $s(\Z\backslash P)$ is a finite totally
  ordered set and therefore there is a unique order preserving
  bijection $t:[n]\rightarrow s(\Z\backslash P)$. Note that
  $t(0)=0$. We can now construct maps $f:[[n]]\rightarrow P$ and $g:P
  \rightarrow [[n]]$ as follows: Let $n_1(i)=[i/(n+1)]$ and let
  $n_2(i)=i- (n+1) n_1(i)$. Then we define
  \begin{equation*}
    f(i)= T_1^{n_1(i)} t(n_2(i)).
  \end{equation*}
  Given $x\in P$, there is a unique $m_1(x)\in \Z$ such that
  $T_1^{m_1(x)}(s\pi(x))=x$ and we define
  \begin{equation*}
    g(x)= t^{-1} s\pi(x)+ (n+1)\cdot m_1(x).
  \end{equation*}
  It is easy to see that $f$ and $g$ are order preserving, commute
  with the shifts and inverse to each other, which completes the proof.
\end{proof}

As for simplicial sets, we may deduce from this a corollary regarding
product of realizations. Note that the product in the category of spaces
with circle action is given by the product of the underlying spaces
together with the diagonal circle action.
\begin{prop}
  If $P$ and $Q$ are compact ppsets, then $||\mathcal{C}(P)\times
  \mathcal{C}(Q)||\isom ||\mathcal{C}(P)||\times 
  ||\mathcal{C}(Q)||$. 
\end{prop}

\end{document}